\documentclass{ifacconf}

\usepackage{xcolor}
\usepackage{graphicx}      
\usepackage{natbib}        
\usepackage{amsfonts}
\usepackage{amsmath}
\usepackage{xcolor}
\usepackage{algorithm}
\usepackage{mathtools}
\usepackage{algpseudocode}
\usepackage{subcaption}
\usepackage{changes}
\definechangesauthor[name=Estefania Loayza, color= orange]{ELR}

\definechangesauthor[name=Zhengang Zhong, color= blue ]{ZZG}

\usepackage{amsmath,tikz}
\usetikzlibrary{matrix, fit}

\newcommand{\R}{\mathbb{R}}

\newcommand{\z}{\mathbf{z}}

\newcommand{\init}{\textup{init}}
\newcommand{\CBOmin}[3]{(#1)^{#2}_{#3}}
\newcommand{\CBOmax}[3]{\left(#1\right)^{#2}_{#3}}
\begin{document}
\begin{frontmatter}

\title{Multi-level Optimal Control with Neural Surrogate Models\thanksref{footnoteinfo}} 

\thanks[footnoteinfo]{This research was supported by the UK Engineering and Physical Sciences Research Council (EPSRC) grant EP/T024429/1 and a Discovery Grant from NSERC (Canada).}

\author[First]{Dante Kalise} 
\author[First]{Estefania Loayza-Romero}
\author[Third]{Kirsten A. Morris} 
\author[Fourth]{Zhengang Zhong}

\address[First]{Department of Mathematics, Imperial College London (e-mail: \{d.kalise-balza,k.loayza-romero\}@imperial.ac.uk)}
\address[Third]{Department of Applied Mathematics, University of Waterloo (e-mail: kmorris@uwaterloo.ca)}
\address[Fourth]{Centre for Process Systems Engineering, Imperial College London (e-mail: z.zhong20@imperial.ac.uk)}

\begin{abstract} 
Optimal actuator and control design is studied as a multi-level optimisation problem, where the actuator design is evaluated based on the performance of the associated optimal closed loop. The evaluation of the optimal closed loop for a given actuator realisation is a computationally demanding task, for which the use of a neural network surrogate is proposed. The use of neural network surrogates to replace the lower level of the optimisation hierarchy enables the use of fast gradient-based and gradient-free consensus-based optimisation methods to determine the optimal actuator design. The effectiveness of the proposed surrogate models and optimisation methods is assessed in a test related to optimal actuator location for heat control.
\end{abstract}


\end{frontmatter}

\section{Introduction}

In engineering, an actuator is a device that materialises the control action within a physical system. Actuators can be mechanical, electrical, hydraulic, or magnetic \citep{kalise2018optimal}. The design of actuators (and sensors) is a fundamental engineering challenge arising, for instance, in vibration control through the design of piezoelectric actuators \citep{peng2005actuator}, and in active noise cancellation in automobiles and aircraft \citep{morris1998noise}.

From a mathematical viewpoint, actuator design precedes the standard control design paradigm where, for a given control system, a control signal is synthesized to achieve a desired performance. Actuator design is concerned with the specification of a control-to-state map which determines the capabilities of the control signal. In this work, we follow an optimisation-based approach for actuator design and control synthesis. We develop a multi-level optimisation pipeline where the actuator design is at highest level of the hierarchy, and is optimised according to the performance of the associated optimal closed loop, which corresponds to the lower level problem. 

We are interested in optimal actuator design problems where the underlying system dynamics are governed by partial differential equations (PDEs). This has been extensively studied in \cite{morris2010linear} for linear quadratic control problems, in \cite{morris2010using, kasinathan2013h} for $H_2$ and $H_{\infty}$ controller design objectives, in \cite{privat} following a spectral approach, and in \cite{kalise2018optimal} in the framework of shape/topology optimisation. In our multi-level approach, this translates into a computationally demanding lower level problem, requiring the solution -in the linear quadratic case- of a large-scale Algebraic Riccati Equation (ARE), which renders evaluations of the higher level problem prohibitively expensive. 

We propose the use of neural network surrogates in our multi-level optimal design/control pipeline to alleviate the computational burden associated to the lower level optimal control problem. The solution of this lower level problem is characterized by a parameter-dependent value function, where the parametric dependence is related to the higher level optimal actuator design. The evaluation of the value function and its gradient, which are intensively used for the parametric optimisation, is approximated by the use of a neural network trained by supervised learning as proposed in in \cite{wei,sara,AKK}. The use of a surrogate for the value function leads to a significant acceleration of the evaluation time for the lower level problem, enabling the use of both gradient-based and gradient-free methods for the solution of the higher level optimisation problems.

The rest of the paper is structured as follows. In Sec.~\ref{sec:LQR}, we introduce the multi-level optimisation framework for optimal actuator and control design.
In Sec.~\ref{sec:NNValueFunction}, we develop neural network surrogates to approximate the value function associated to the lower level of our hierarchy. We propose the construction of both unstructured and structured surrogates, and evaluate their performance in a prototypical example. Having built a suitable surrogate, Sec.~\ref{sec:algorithms} is devoted to the presentation of optimisation algorithms to solve the max-min problem arising in the higher level of our hierarchy. We present a projected gradient descent ascent and a consensus-based method for saddle point problems and we illustrate their performance in a problem related to optimal actuator location for a thermal system.

\section{A multi-level optimisation framework for optimal actuator/control design}
\label{sec:LQR}
We consider parameter-dependent linear dynamical systems described by
\begin{equation}
\label{eq:system_general}
\frac{\d z}{\d t}=A z(t)+B(r) u(t), \quad z(0)=z_0
\end{equation}
where $A\in \mathbb{R}^{n \times n}$, and $B(r)\in\mathbb{R}^{n\times m}$ is matrix-valued function depending on a parameter $r\in\mathbb{R}^m$. Such systems naturally arise after semi-discretization in space of systems governed by PDEs, where $r$ parametrizes the location of $m$ actuators. For the sake of simplicity, we assume that the control space is of the same dimension of the parameter space. However, the methodology developed in this work can be seamlessly adapted to work over spaces of different dimensions. We assume that each parameter coordinate can be varied over some compact set $\Omega \subset$ $\mathbb{R}$.

The linear quadratic controller design aims at finding a minimising control $u(t)\in \mathbb{R}^m$ to the cost functional  
\begin{equation}
    \label{eq:obj_LQR}
J\left(u; z_0,r\right)=\int_0^{\infty} z(t)^\top Q z(t) +  u(t)^\top R u(t),
\end{equation}
where $Q \in \mathbb{R}^{n\times n}$, $Q\succeq 0$, $R \in \mathbb{R}^{m\times m}$, $R\succ 0$,  and $z(t)\in \mathbb{R}^n$ is determined by the dynamics \eqref{eq:system_general}.

For each parametric realisation $r$ we assume that the pairs $\left(A, B(r)\right)$ are stabilisable. 
The optimal cost-to-go or value function for a given initial condition $z_0$ and parameter $r$ is 
\begin{equation}
\label{eq:value_function}
V\left(z_0,r\right)\coloneqq \inf_{u \in \mathcal{U}} J(u;z_0,r) =  z_0^\top \Pi(r) z_0,
\end{equation}
 with $\mathcal{U} = L^2\left((0,+\infty);\mathbb{R}^m\right)$,

where $\Pi(r)$ solves the parameter-dependent Algebraic Riccati Equation (ARE)
\begin{equation}
    \label{eq:general_parameterized_ARE}
A^\top \Pi(r)+\Pi(r) A-\Pi(r) B(r) R^{-1} B(r)^\top \Pi(r)+ Q=0.
\end{equation}
We can define a parametric optimisation problem, that is,
\begin{equation}
\min_{r \in \Omega^m} \, V(z_0,r) =  \min_{r \in \Omega^m} z_0^\top \Pi(r) z_0,
\end{equation} 
however, in general, the initial condition $z_0$ is not fixed.  For example, we can optimise $r$ according to the worst possible initial condition, that is,
\begin{equation}
    \label{eq:general_max-min_opt_prob}
\max _{\substack{z_0 \in \mathbb{R}^n \\\left\|z_0\right\|=1}} \min _{u \in \mathcal{U}} J\left(u; z_0, r\right)=\max _{\substack{z_0 \in \mathbb{R}^n \\ \left\|z_0\right\|=1}} z_0^\top \Pi(r) z_0,
\end{equation}
where the last expression corresponds to the Rayleigh quotient, i.e., $\|\Pi(r)\|$.
We refer the reader to \cite{morris2010linear,edalatzadeh2021optimal,kalise2018optimal} for more details on how to solve this problem. 

Overall, the simultaneous optmisation of the control signal and the parameter, together with an optimality-based characterisation of the space of initial conditions of interest, induces a hierarchy of costs which are cast as a multi-level optimisation problem

\begin{equation}
\label{eq:max-min_fin_dim}
    \max _{z_0 \in \mathbb{R}^n} \min_{r\in \Omega^m} \min _{u \in \mathcal{U}} J(u;z_0,r) =  \max _{z_0 \in \mathbb{R}^n} \min_{r\in \Omega^m} V(z_0,r).
\end{equation}

To eliminate the burden of computing the solution of the ARE ~\eqref{eq:general_parameterized_ARE} at every evaluation of $V(z_0,r)$ within an algorithm for the outer max-min problem, we propose to build a closed-form surrogate for $V(z_0,r)$ using neural networks.

\section{Surrogates for the value function with Neural Networks}
\label{sec:NNValueFunction}

In a nutshell, deep learning is about realising complex tasks such as speech, or image recognition, and language translation, among others, by means of parameterised functions called neural networks (NNs).
A neural network consists of neurons that are ordered into layers. We have three different kinds of layers: input, hidden and output layers. 
Such a network is called multilayer perceptron (MLP)\citep{svozil1997introduction}. 
Here, we are interested in building neural network surrogates $V_{\theta}(z_{0}, r): \mathbb{R}^{n + m} \rightarrow \mathbb{R}$ for the value function $V(z_{0}, r)$, described in~\eqref{eq:value_function}.

In what follows, all the neural networks considered are single hidden layer feedforward NN with \texttt{ReLu} and/or \texttt{softplus} as activation functions.

\subsection{Unstructured Surrogates}
\label{subsec:directSurrogate}

To train the neural network, we collect $N_z \times N_r$ joint pairs for the initial state $z_{0}$ and the design parameter $r$ from a uniform sampling of $\| z_0\| \le 1$ and $\Omega^m$, respectively. To generate output data, we calculate $\hat{y}:= z_{0}^\top \Pi(r) z_{0}$ for each pair of initial state and design parameter. To train the neural network parameters $\theta$, we consider the mean squared error loss function  
\begin{equation}
\label{eq:lossFunctionDirect}
    \mathcal{L}(\theta)=\frac{1}{\mathrm{N_z \times N_r}} \sum_{\mathrm{i}=1}^{\mathrm{Nz \times N_r}}\left(V_{\theta}(z_{0}^{(i)}, r^{(i)})-\hat{y}^{(i)}\right)^2.
\end{equation} 

It is worth highlighting that from optimal control theory one expects the value function to be non-negative since the solution of the ARE is positive semi-definite. However, with this unstructured approach, this cannot be guaranteed. In the following, we propose a structured surrogate model by approximating the matrix $\Pi(r)$ instead.

\subsection{Structured Surrogates}
\label{subsec:NNPiMatrix}
As the value function \eqref{eq:value_function} is a quadratic function of the initial state and requires the solution of the ARE \eqref{eq:general_parameterized_ARE}, we consider alternatively to construct a NN directly approximating $\Pi(r)$, which is a positive semi-definite matrix, thus guaranteeing the non-negativeness of the value function.

We define a NN surrogate model for the solution as $\Pi(r) \approx \Pi_{\theta} := L_{\theta}(r)L_{\theta}^{\top}(r)$, where $L_{\theta}(r)$ aims at learning the unique Cholesky decomposition of symmetric positive definite matrices. 
Let $NN_{ij}$ denote a MLP as defined in
\citep{svozil1997introduction}
mapping $r$ to a scalar value in $ij$-th element of the matrix $L_{\theta}$. We define the matrix-valued surrogate model as 
\begin{equation}
    f_{ij} =
    \begin{cases}
       \varphi(NN_{ij}) + \varepsilon, & \text{if } i = j, \\ 
        NN_{ij}, & \hbox{otherwise.} \\
    \end{cases}
\end{equation}

where $\phi$ is an activation function. As it is known that $L_{\theta}(r)$ has positive diagonal entries, we define the neural network surrogate as

\begin{center}
\begin{tikzpicture}[baseline=(current bounding box.center)]
\matrix (m) [matrix of math nodes, nodes in empty cells,right delimiter={]},left delimiter={ [ } ]{
f_{11}(r)  & & & \\
 & & & \\
 & & &   \\
f_{n1}(r) & &  & f_{nn}(r) \\
} ;
\node[fit=(m-1-3)(m-3-4)]{\hspace*{20pt}\Large{0}};
\draw[loosely dotted] (m-1-1) -- (m-4-4);
\draw[ dotted] (m-4-1) -- (m-4-4);
\draw[ dotted] (m-1-1) -- (m-4-1);
\end{tikzpicture}
\end{center}
where off-diagonal elements take a two-layer structure mapping the design parameter to a scalar. 

For the purpose of training, we generate $N$ design parameter samples from $\Omega^m$ and thereby the corresponding Riccati solutions $\Pi(r_{i})$. The loss function is hence defined as 
\[
\mathcal{L}(\theta) = \sum_{i = 1}^{N} \|\Pi_{\theta}(r_i) - \Pi(r_{i})\|^2_{\textup{fro}},
\]
where $\|\cdot\|_{\textup{fro}}$ is the Frobenius norm. 

\subsection{Numerical verification}\label{nv}
We consider a prototypical example related to optimal actuator location problem for a heating system, similarly as in \cite{kalise2018optimal}. After semi-discretization in space, the finite dimensional linear dynamics read
\begin{align}
\nonumber 
\frac{\d z}{\d t} & = A z(t) + B(r) u(t), & t\in [0,+\infty[ ,\\
z(0)&=z_0, &
\label{eq:heat_equation}
\end{align}
where $A \coloneqq KM^{-1}$, with
\begin{align*}
    K_{ij} &= \int_0^\pi \phi_i'(x) \phi_j'(x) \ \d x &
    M_{ij} &= \int_0^\pi \phi_i(x) \phi_j(x) \ \d x \\
    B(r)_{i\ell}   &= \int_{r_\ell-\delta\pi}^{r_\ell + \delta\pi} \phi_i(x) \ \d x \ ,
\end{align*}
with the functions $\phi_i(x) = \sin(ix)$ on $[0,\pi]$, for all $i,j= 1\ldots, n$ and $\ell = 1\ldots, m$. The matrices appearing in~\eqref{eq:obj_LQR} are given by $Q = M$ and $R = \textup{id}_{m\times m}$. Here, the parameter $r$ represents the location of $m$ patches of length $2\delta$ inside the domain, each one associated to a scalar control signal.

In the following, we demonstrate the performance of unstructured and structured surrogates $V_{\theta}(z_{0},r)$ for a problem with $n = 3,\delta=0.005,$ and $ m = 1$, i.e., a single actuator. 
Consider $\hat{z}_0:= \prod_{j \in 1\dots n}\mathcal{Z}_j$, where $\mathcal{Z}_j := \{- 1 + (i-1)/3.5 \mid i \in (1, \dots, 8)\}$ and $\hat{r}:= \{i\pi/100 \mid i \in (1, \dots, 99)\}$. We calculate the value $z_0^\top \Pi_{n}(r)z_0 $ from $\{(z_0, r) \mid z_0 \in \hat{z}_0, r \in \hat{r} \}$, thereby $51200$ input and output data for training the neural network. 
We construct the surrogate neural network with one hidden layer with $128$ neurons and \texttt{ReLu} as the activation function. 
The training optimisation is solved with the \texttt{Adam} solver with the learning rate $1e-3$ and $223000$ iterations. 
In this work, all NN surrogates are built by using \texttt{PyTorch}. 
 
\begin{figure}
\begin{center}
\includegraphics[width=8.4cm]{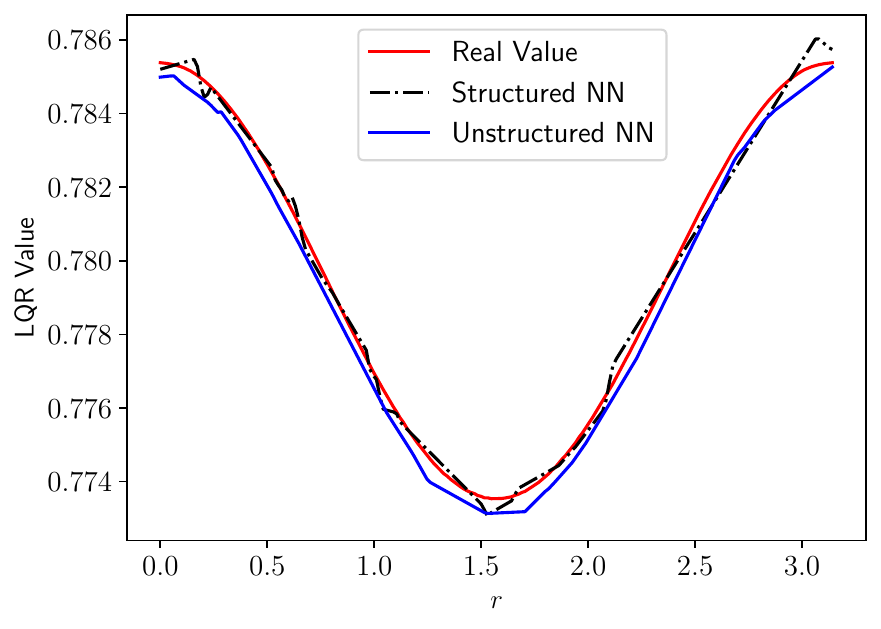} 
\caption{True $\max_{\|z_{0}\| = 1}V(z_{0}, r)$, structured and unstructured surrogates $\max_{\|z_{0}\| = 1}V_{\theta}(z_{0}, r)$ for $m=1$ and $n = 3$. Both surrogates effectively approximate the value function without violating the non-negativeness constraint, the structured surrogate generates a more accurate solution.}
\label{fig:5dsys_true_worst_value}
\end{center}
\end{figure}

For the structured surrogate we compute $\Pi(r)$ at $r \in \hat{r}:= \{i\pi/120 \mid i \in (1, \dots, 119)\}$, thereby $120$ solutions of the Riccati equation are applied for training the neural network. 
Each element of the surrogate matrix takes MLP with one hidden layer with $128$ neurons. 
It uses \texttt{ReLu} as the activation function and $\varphi$ as \texttt{softplus} (function guaranteeing the positive definiteness) with $\varepsilon \approx 1.192e-7 $, i.e. the default smallest representable number for \texttt{float32} in \texttt{PyTorch}.
The training optimisation is solved with the \texttt{Adam} solver with the learning rate $1e-3$ and $6000$ iterations. 
In Fig.~\ref{fig:5dsys_true_worst_value}, we show the values of the real value function and the approximated unstructured and structured NN surrogates. It can be seen that while both surrogates effectively approximate the value function without violating the non-negativeness constraint, the structured surrogate generates a more accurate solution.

\begin{figure}
\begin{center}
\includegraphics[width=8.4cm]{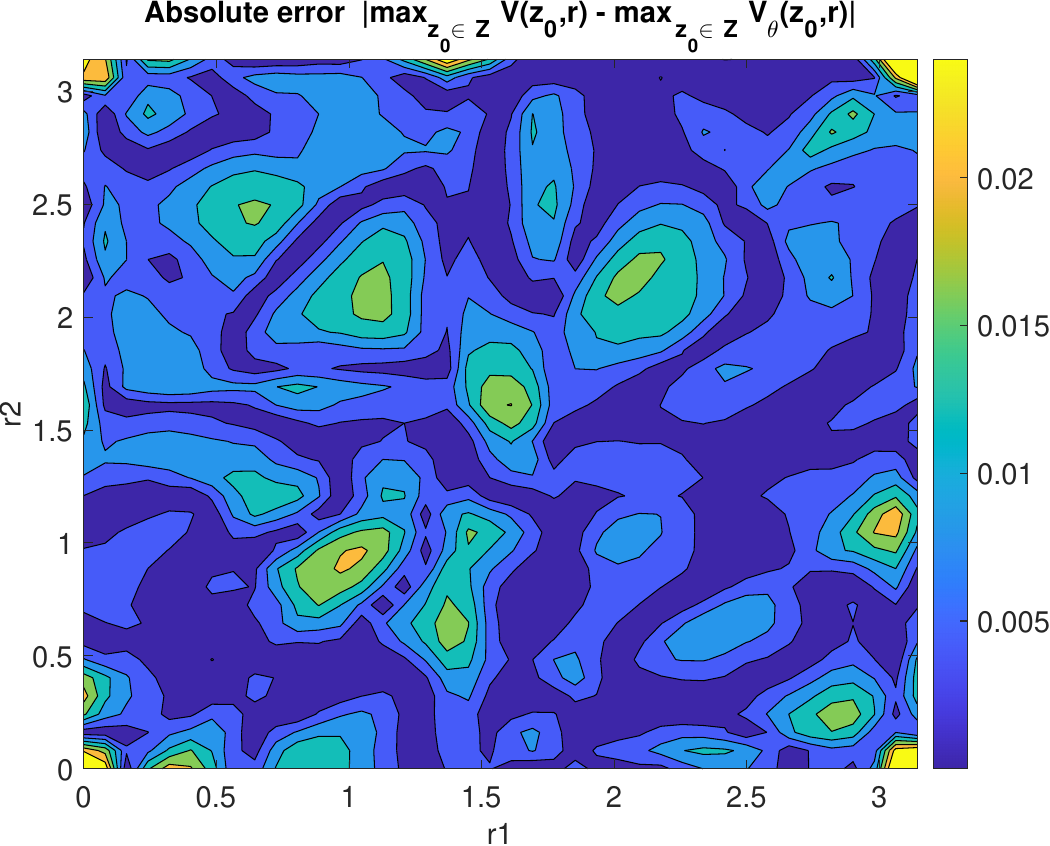}    
\caption{The heat map of absolute training error using the structured NN for $n =10$ and $m = 2$, with the surrogate constructed as indicated in Subsec.~\ref{subsec:NNPiMatrix}. }
\label{fig:10dsys_error}
\end{center}
\end{figure}

We further consider the performance of the structured surrogate $V_{\theta}(z_{0},r)$ on an approximation of \eqref{eq:heat_equation} with $n =10, m = 2$. We calculate $\Pi(r)$ at $r := (r_1, r_2) \in \{ (r_1, r_2) \mid r_1, r_2 \in \hat{r}\}$, where $\hat{r}:= \{(i-1)\pi/19 \mid i \in (1, \dots, 20)\}$, thereby $400$ solutions of the Riccati equation are applied for training the neural network. Fig.~\ref{fig:10dsys_error} shows the heat map corresponding to the absolute error of $| \max_{\|z_0\|=1}V(z_{0},r) - \max_{\|z_0\|=1} V_{\theta}(z_{0},r) |$, from where it can be observed that the surrogate yields an accurate approximation of the true value function.

The training process of the unstructured surrogate is limited by computational challenges, as it requires significantly more training data -- $(N_x)^{n}$ times more -- compared to the structured surrogate for the same actuator location grids. In the following we restrict ourselves to the construction of structured surrogates.

\section{Surrogate Optimisation}

\label{sec:algorithms}

The true value function $V(z_0,r)$ and its gradient are expensive to evaluate as they require the solution of an ARE and additional sensitivity relations in the case of the gradient. Instead, the NN surrogate $V_{\theta}(z_0,r)$ and its gradient have closed-form expressions that can be readily evaluated using automatic differentiation. In the following, we use the trained surrogate to develop gradient-based and gradient-free algorithms for the solution of the max-min optimisation problem 
\begin{equation}
    \label{eq:max-min_problem}
    \max_{z_{0} \in Z}\min_{r \in \Omega^m} V_{\theta}(z_{0},r)
\end{equation}
with $Z = \{z \in \R^n \colon \lVert z \rVert= 1\}$ and $\Omega^m = [0,\pi]^m$.

\subsection{Projected gradient descent ascent (PGDA)}

One natural candidate for solving problem \eqref{eq:max-min_problem} is the gradient descent-ascent (GDA) algorithm, which, at each iteration, performs gradient descent over the variable $r$ with the step-size $\eta_{\mathrm{r}}$ and gradient ascent over the variable $z_{0}$ with the step-size $\eta_{z_{0}}$. Moreover, thanks to the way the structured surrogate is constructed, we can guarantee that the value function $V_\theta$ is convex w.r.t. $z_0$; thus the non-convex constraint $Z= \{z \in \R^n \colon \lVert z \rVert = 1\}$ can be safely replaced by the convex set $\bar{Z} = \{z \in \R^n \colon \lVert z \rVert \leq  1\}$ without altering the solutions of the problem. 

Since the related projection operators onto $Z$ and $\Omega^m$ can be computed with little effort, we consider a projected version of the gradient descent ascent (PGDA).
We recall, the projection operators are given by $\operatorname{Proj}_{\bar{Z}}(z) = z/||z||$, if $||z||>1$ and $\operatorname{Proj}_{\bar{Z}}(z) = z$ if $||z||\leq 1$. On the other hand, 
$\operatorname{Proj}_{\Omega^m}(r) = \max\{\min\{r_{ub},r\},r_{lb}\}$, where $r_{ub}\in \mathbb{R}^m$ has all its entries equal to $\pi$ (upper bound of the box constrains) and $r_{lb}\in \mathbb{R}^m$ is the zero vector (lower bound).  
Here the operations $\min$ and $\max$ are acting in component-wise. The method is described in  Alg.~\ref{alg:PGDA}. 

\begin{algorithm}
\caption{Projected gradient descent-ascent}
\label{alg:PGDA}
\begin{algorithmic}[1] 
  \Procedure{PGDA}{$z_0^{\init},r^{\init}, \eta_{\z_0}, \eta_{r}, Z, \Omega^m$}       
  \Require Objective $V_{\theta}$, initial guess $z_0^{\init},r^{\init}$, step-size $ \eta_{z_0}, \eta_{r}$, and feasible set $Z, \Omega^m$  
  \State $k \gets 0$
  \State $z_0^{k} \gets z_0^{\init}$ and $r^{k} \gets r^{\init}$
  \While{$k \leq K$}  
      \State $z_0^{k+1}=\operatorname{Proj}_{\bar{Z}}\left[z_0^k + \eta_{z_{0}} \nabla_{z_0} V_\theta\left(z_0^k,r^k\right)\right]$
      \State $r^{k+1}=\operatorname{Proj}_{\Omega^m}\left[r^k-\eta_r \nabla_{r} V_\theta\left(z_0^k,r^k\right)\right]$
      \State $k \gets k+1$
    \EndWhile  
    \State \textbf{return} $z_0^{k}, r^{k}$
  \EndProcedure
\end{algorithmic}
\end{algorithm}

Even when considering linear dynamics, using NN surrogates to approximate the value function may result on a lack of convexity-concavity. 
Moreover, it is well known that PGDA fails to converge to a global Nash equilibrium, see e.g.,~\citep{huang2022consensus}. 
Therefore, in this study, we aim to overcome this limitation by employing a consensus-based algorithm, as proposed by \citep{huang2022consensus}, to address the max-min surrogate optimisation problem \eqref{eq:max-min_problem}. 

\subsection{Consensus-based method for saddle points}

A general consensus-based algorithm (CBO) is an agent-based global optimisation algorithm, which mimics interacting agents communicating over a weighted mean \citep{totzeck2021trends}. CBO methods use a finite number of agents, which are formally stochastic processes, to explore the domain and to form a global consensus about the minimiser as time passes. The dynamics of the agents are governed by two competing terms. A drift term that drags the agents towards a momentaneous consensus point
and a diffusion term that randomly moves agents according to a scaled Brownian motion, featuring the exploration of the energy landscape of the cost.
CBO compared to other agent-based methods, such as Particle Swarm Optimisation (PSO), allows rigorous analysis of convergence when it comes to the mean-field limit. For $N$ particles, the communication is of order $\mathcal{O}(N)$. In the following, we summarize the consensus-based method for saddle points proposed in \cite{huang2022consensus}.

In the max-min setting, two sets of agents corresponding to the minimiser $\{\CBOmin{r}{i}{}\}_{i = 1}^{N_1}$ and maximiser $\{\CBOmax{z_0}{i}{}\}_{i = 1}^{N_2}$ are considered, respectively. To achieve consensus about the equilibrium point of the interested max-min value function, the agents interact through a system of stochastic differential equations (SDEs) of the form 
\begin{align*}
\d \CBOmin{r}{i}{t} & =-\lambda_1\left(\CBOmin{r}{i}{t}-x_\alpha^{z_0}\right) \d t+\sigma_1 D\left(\CBOmin{r}{i}{t}-x_\alpha^{z_0}\right) \d W_t^{r, i},\\
\d \CBOmax{z_0}{i}{t} & =-\lambda_2\left(\CBOmax{z_0}{i}{t}-y_\beta^r\right) \d t+\sigma_2 D\left(\CBOmax{z_0}{i}{t}-y_\beta^r\right) \d W_t^{z_0, i}.
\end{align*}
To improve readability we have omitted the dependence of the variables $x_\alpha^{z_0}$ and $y_\beta^r$ on $\left(\widehat{\rho}_{r, t}^{N_1}\right)$ and $\left(\widehat{\rho}_{z_0, t}^{N_2}\right)$, respectively. Moreover, we have 
$\widehat{\rho}_{r, t}^{N_1}=\frac{1}{N_1} \sum_{i=1}^{N_1} \delta_{\CBOmin{r}{i}{t}}$, and  $\widehat{\rho}_{z_0, t}^{N_2}=\frac{1}{N_2} \sum_{i=1}^{N_2} \delta_{\CBOmax{z}{i}{t}}$, and $D(\cdot)$ stands for $\operatorname{diag}(\cdot) $. 
 
The initial conditions $\CBOmin{r}{i}{0} \sim \rho_{r, 0} \in \mathcal{P}\left(\mathbb{R}^{d_1}\right)$ for $i=1, \ldots, N_1$ and $\CBOmax{z_0}{i}{0} \sim \rho_{z_0, 0} \in \mathcal{P}\left(\mathbb{R}^{d_2}\right)$ for $i=1, \ldots, N_2$. 
The independent standard Brownian motions in $\mathbb{R}^{n}$ and $\mathbb{R}^{m}$ are denoted by $\left(W_t^{r, i}\right)$ for $t\geq 0$ and $i=1,\ldots, N_1$ and $\left(W_t^{z_0, i}\right)$, for $t\geq 0$ and $i=1,\ldots, N_2$.

Moreover, $\widehat{\rho}_{r, t}^{N_1}$ and $\widehat{\rho}_{z_0, t}^{N_2}$ denote the empirical measures of the particles' $x$- and $y$-positions, respectively. 
The consensus point $\left(x_\alpha^{z_0}\left(\widehat{\rho}_{r, t}^{N_1}\right), y_\beta^r\left(\widehat{\rho}_{z_0, t}^{N_2}\right)\right)$ is defined by 
\begin{subequations}
\begin{align}
x_\alpha^{z_0} & =\int x \frac{\omega_\alpha\left(x, \int y d \widehat{\rho}_{z_0, t}^{N_2}(y)\right)}{\left\|\omega_\alpha\left(\cdot, \int y d \widehat{\rho}_{z_0, t}^{N_2}(y)\right)\right\|_{L_1\left(\widehat{\rho}_{r, t}^{N_1}\right)}} \d \widehat{\rho}_{r, t}^{N_1}(x), \label{eq:consensus_sub1}\\
y_\beta^r &=\int y \frac{\omega_{-\beta}\left(\int x d \widehat{\rho}_{r, t}^{N_1}(x), y\right)}{\left\|\omega_{-\beta}\left(\int x d \widehat{\rho}_{r, t}^{N_1}(x), \cdot\right)\right\|_{L_1\left(\widehat{\rho}_{z_0,t}^{N_2}\right)}} \d \widehat{\rho}_{z_0, t}^{N_2}(y), \label{eq:consensus_sub2}
\end{align}
\end{subequations}
where $\omega_\alpha(x, y):=\exp (-\alpha V_{\theta}(y, x))$ and $\omega_{-\beta}(x, y):=\exp (\beta V_{\theta}(y, x))$.

The agent simulations are realised using the standard Euler-Maruyama time discretisation with the step size $\Delta t$, i.e.
\begin{subequations}
\begin{align}
\CBOmin{r}{i}{k+1}&=\CBOmin{r}{i}{k}-\lambda_1 \Delta t\nu(r)+\sigma_1 D\left(\nu(r)\right) W_k^{r, i}, \label{eq:update_discrete_sub1} \\
\CBOmax{z_0}{i}{k+1}&=\CBOmax{z_0}{i}{k}-\lambda_2 \Delta t\nu(z_0)+\sigma_2 D\left(\nu(z_0)\right) W_k^{z_0, i},
\label{eq:update_discrete_sub2}
\end{align}
\end{subequations}
where $\nu(r) = \CBOmin{r}{i}{k}-x_\alpha^{z_0}$ and $\nu(z_0) = \CBOmax{z_0}{i}{k}-y_\beta^r$. Notice that we use the same notation for the continuous and discretised versions of the system describing the interaction of the agents.

The algorithm applied in this work sequentially updates the discredited states and consensus points via \eqref{eq:update_discrete_sub1}, \eqref{eq:update_discrete_sub2}, \eqref{eq:consensus_sub1}, \eqref{eq:consensus_sub2} until it converges as proposed by \citep[Algorithm 1]{huang2022consensus}.

To penalise iterates leaving the set $Z \times \Omega^m$, we use a modified objective function, i.e. penalty with a positive value for the minimiser and with a negative value for the maximiser in the objective 
$V := V_{\theta} + \mu(-|\|z_0\| - 1| + \sum_{l}([r - r_{ub}]_{+}^{l} + [r_{lb} - r]_{+}^{l})),$
where $\mu$ is a scalar penalizing the constraint violation, $[\cdot]_{+}^{l}$ is the $l$-th component of the positive part of a vector and $r_{lb}$, $r_{ub}$ are low and upper bound of the box constraints of $r$, respectively.

\subsection{Numerical tests}

\begin{figure}
\begin{center}
\includegraphics[width=8.4cm]{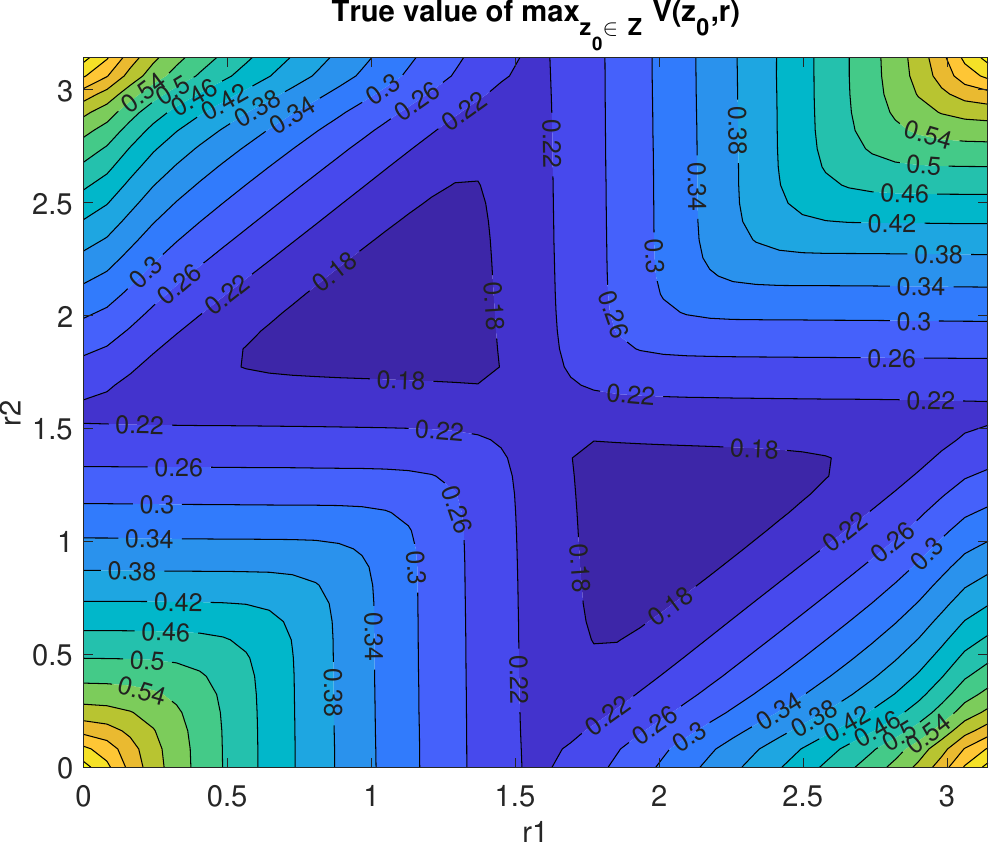}    
\caption{The heat map corresponding to the real worst-case value function.} 
\label{fig:10dsys_true_value}
\end{center}
\end{figure}

\begin{figure}
\begin{center}
\includegraphics[width=8.4cm]{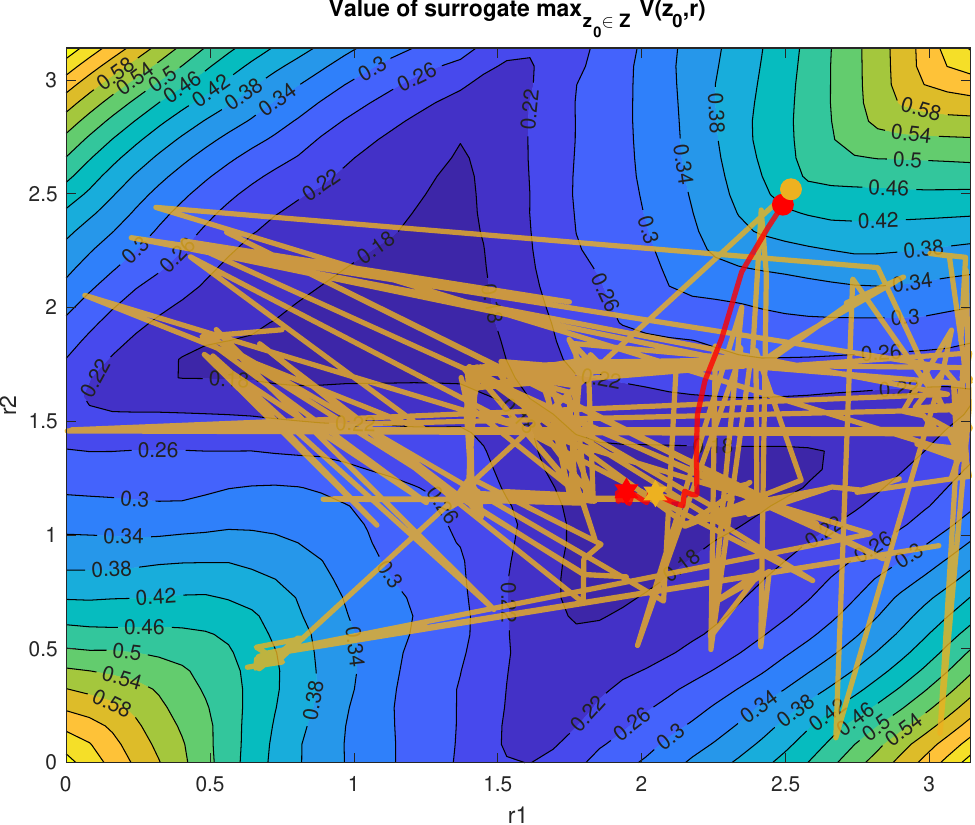}    
\caption{The heat map corresponding to the NN surrogate worst-case value function. Red: solution trajectories of PGDA, Magenta: consensus point trajectories of CBO-SP. Initial points are depicted with a circle and end points with a star.}
\label{fig:10dsys_surrogate_value}
\end{center}
\end{figure}

We return to the example presented in Section \ref{nv} related to optimal actuator location. We set $n = 10, m = 2$ and we solve the surrogate max-min problem \eqref{eq:max-min_problem}. The heat maps corresponding to the real and surrogate worst-case value functions are shown in Figs.~\ref{fig:10dsys_true_value} and~\ref{fig:10dsys_surrogate_value}, respectively.

For PGDA we select $K = 2000$, $\eta_{r} = 3\times 10^{-4}$, $\eta_{\z_0} = 10^{-3}$ and initial guess $z_0^{\init} = [0.5, \cdots, 0.5], r_\ell^{\init} = [2.5, 2.5]$.
For CBO-SP, we set $K = 2000$, $\alpha = \beta = 10^{15}, \lambda_1 = 2, \lambda_2 = 0.1, \sigma_1 = \sigma_2 = 2$, $\mu = 10000$ and initial law corresponding to normal distributions $\rho_{r,0} = \mathcal{N}(r_0^{\init}, diag(1.5,1.5))$, $\rho_{z_0,0} = \mathcal{N}(z_0^{\init}$, $diag(1.5, \cdots, 1.5)), N_1 = N_2 = 300.$

The solutions obtained by PGDA and CBO-SP after $K$ iterations are $r_{CBO-SP} = [2.0472, 1.1830]$ and $r_{PDGA} = [1.9465, 1.1927]$, respectively, are displayed Fig.~\ref{fig:10dsys_surrogate_value}. Both algorithms converge to similarly accurate solutions, distributing the actuators symmetrically in $[0,\pi]$. From the surrogate value function it can also be observed that actuators can be interchanged, as expected.

We validate the closed-loop performance of the optimal parameter and associated optimal control signals $r_{CBO-SP}$. For reference, we consider a sub-optimal solution where both actuators are placed at the origin $r = [0, 0]^{\top}$.

\begin{figure}
\begin{center}
\includegraphics[width=9cm]{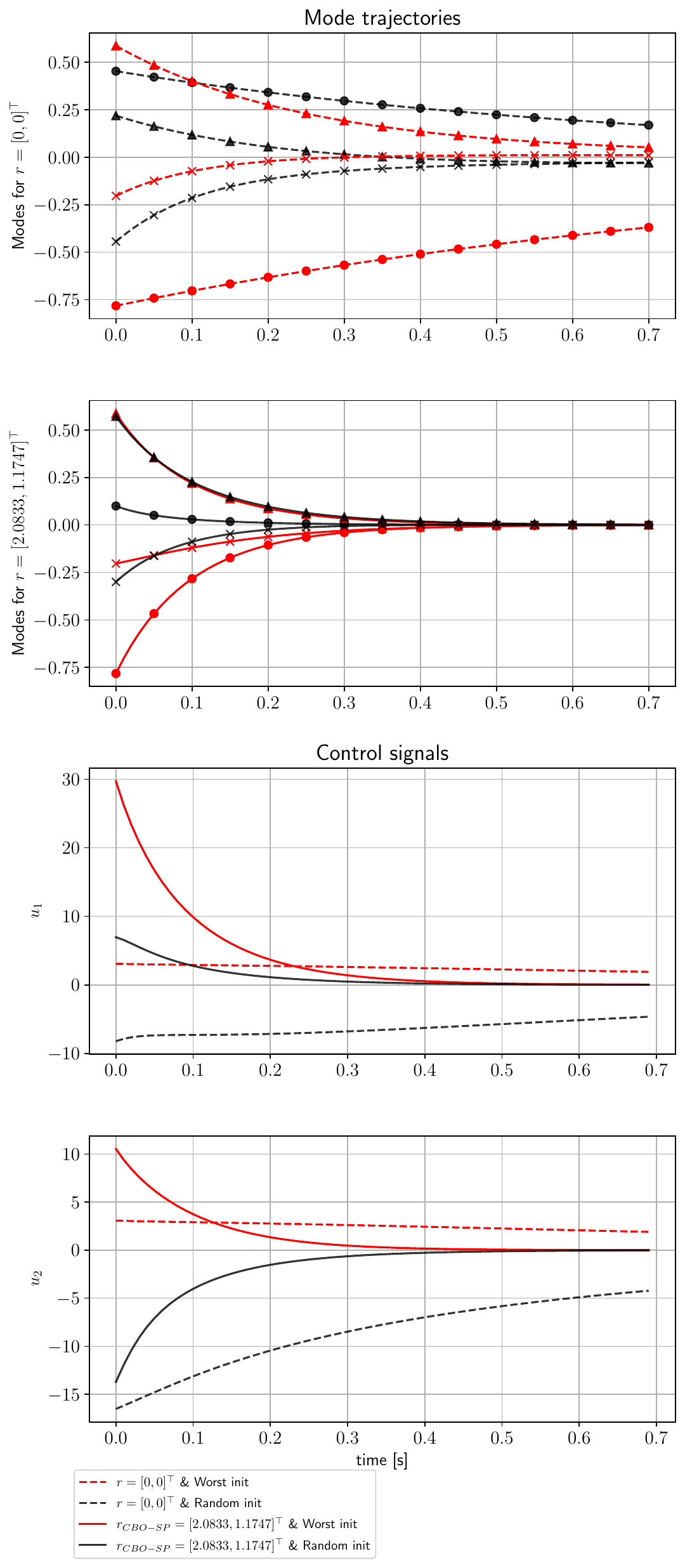}    
\caption{Closed-loop performance of the optimised actuators at $r_{CBO-SP} = [2.0472, 1.1830]$ and non-optimised ones placed at $r = [0, 0]^{\top}$. The first three modes are demonstrated with circle, triangle and cross, respectively.}
\label{fig:3d_plot_opt_design}
\end{center}
\end{figure}

In Fig.~\ref{fig:3d_plot_opt_design}, we observe that the state of the first mode converges significantly slower when the actuators are placed in a sub-optimal position. Additionally, in Fig.~\ref{fig:3d_plot_opt_design} we show that under the worst-case initial condition obtained with CPO-SP, the system states can be controlled to zero within $0.7$ seconds with an optimal actuator position, whereas the system states remain far from zero under a sub-optimal actuator location.

\section{Conclusion and Future Work}
We have shown that multi-level control problems related to optimal actuator/sensor design, where the optimal parameter is linked to an optimal closed-loop performance, can be computationally alleviated by the use of neural network surrogates which. This surrogate can be endowed with additional structures which are relevant to the optimisation, such as non-negativeness or convexity. Moreover, the use of a neural network surrogate with a closed-form expression after training, allows fast evaluation of the low-level value function and its gradient using automatic differentiation, which is considerably cheaper than sensitivity expressions for the value function requiring adjoint calculations. This enables the use of both gradient-based and gradient-free solvers for the higher-level max-min problem.

Our main interest in this subject is the development of an effective pipeline for optimal actuator design and location in nonlinear distributed parameter systems. We note that while some parts of the proposed methodology may seem redundant for the linear-quadratic presented in this work, such as the use of a gradient-free solver for the max-min problem, they will become essential in the nonlinear case, where we can no longer rely on a Riccati-based approach to characterize the value function of the low-level optimal control problem.

\bibliography{ifacconf}             
                                                   
\end{document}